\documentclass{article}	
\usepackage{float}		
\usepackage{times}		
\usepackage{graphicx}		
\usepackage{url}
\usepackage{amsmath}
\usepackage{amssymb}
\usepackage{cite}

\usepackage{booktabs}
\usepackage[table]{xcolor}

\def\EPSFIG[#1]#2#3#4{		
\begin{figure}[H]		
\begin{center}			%
\includegraphics[#1]{#2}	%
\end{center}			%
\caption{#3}			%
\label{#4}			%
\end{figure}			%
}				%

\author{Javier L\'opez Pe\~na\footnote{Department of Mathematics, University College London. E-mail: j.lopezpena@ucl.ac.uk} \ and
Hugo Touchette\footnote{School of Mathematical Sciences. Queen Mary University of London. E-mail: h.touchette@qmul.ac.uk}}
\title{A network theory analysis of football strategies}

\begin{document} 
\twocolumn
\maketitle 
\thispagestyle{empty}


\abstract{We showcase in this paper the use of some tools from network theory to describe the strategy of football teams. Using passing data made available by FIFA during the 2010 World Cup, we construct for each team a weighted and directed network in which nodes correspond to players and arrows to passes. The resulting network or graph provides a direct visual inspection of a team's strategy, from which we can identify play pattern, determine hot-spots on the play and localize potential weaknesses. Using different centrality measures, we can also determine the relative importance of each player in the game, the `popularity' of a player, and the effect of removing players from the game.}

\section{Introduction}
Graphs or networks arise in the study of a variety of problems, ranging from technological and transport issues to social phenomena and biological problems \cite{Lewis2009, Newman2010, Boccaletti2006}. Their prevalence is such that a rich mathematical theory has been developed around them, notably by Euler, in relation to the K\"onigsberg bridge problem, Erd\"os and many others.

In the world of sports, team sports involving passes between players provides one with interesting examples of networks. Our goal in this paper is to show how the mathematical theory of networks can be used to analyze statistical information of team sports and measure the performance of a team and its players. As a proof of concept, we apply our ideas to construct a network analysis of some of the teams participating in the football 2010 World Cup.

Arguably the most popular sport in the world, football (\emph{soccer} for our American readers) has traditionally lagged behind other sports, such as baseball or basketball, in terms of statistical information made available after games. The unique nature of football games, with their constant ball flow and comparatively low scores compared to other sports, makes simple statistics such as assists or number of goals insufficient as measures of team and player performance. Fortunately, it seems that the situation is changing. In recent years, starting with the UEFA 2008 Euro Cup, an unprecedented amount of statistical data has been made public after games. The release of significantly larger amounts of data opens up the way for building new and more detailed analyses of football. Some recent attempts in this direction can be found in \cite{Duch2010,Brillinger2007,Hughes2005,Joseph2006,Luzuriaga2010,Yamamoto2011}.

Here we focus on finding a quantifiable representation of a team's style using network theory. All renowned football teams in history have displayed a recognizable footprint in their game-play, which has always been thought of as something observed by football experts rather than described by game statistics. To reveal this footprint, we use the passing distribution of a team to construct a weighted and oriented network, with nodes corresponding to players and  weighted arrows to the number of successful passes between players. By attaching each node to the tactical positioning of the team, we then obtain an immediate picture of the team's style, which can profitably be used to observe overused and underused areas of the pitch or to detect potential performance problems between certain players. By computing certain network invariants, such as centrality measures, we can also analyze a team's performance as well as the contributions of each of its players. These measures yield, as will be seen, a lot of useful information despite the relatively small size (11 players) of passing networks.

\section{The network of a football team}

\begin{figure*}
  \begin{center}
    \includegraphics[width=0.8\textwidth]{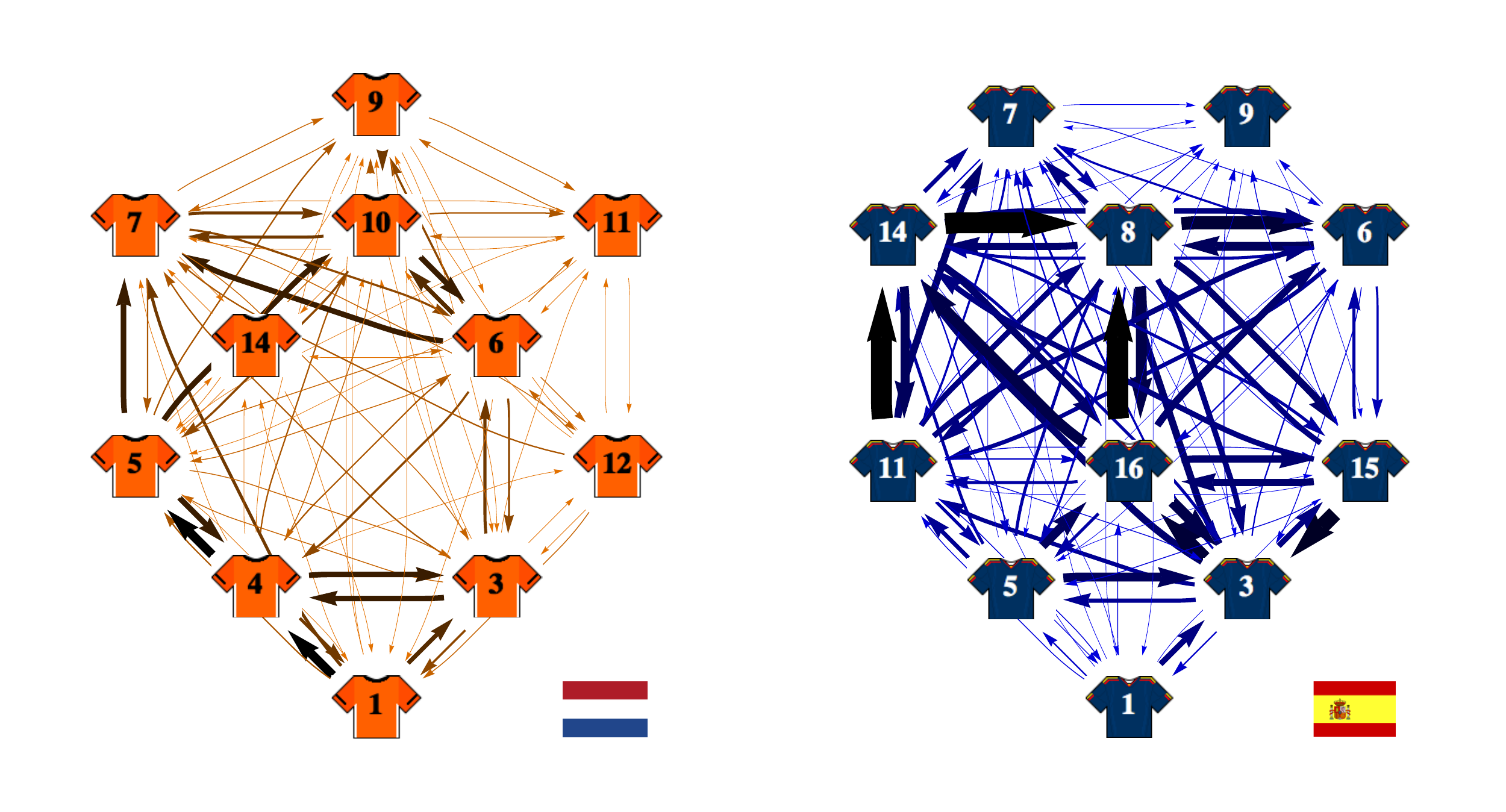}
  \end{center}
  \caption{Passing networks for the Netherlands and Spain drawn before the final game, using the passing data and tactical formations of the semi-finals.}
  \label{fig:spvsned}
\end{figure*}

We define the \emph{passing network} of a team as the network containing the team's players as nodes and connecting arrows between two players weighted by the successful number of passes completed between them. Although networks are, technically speaking, only topological in nature, we use the passing network as a tool for visualizing a team's strategy by fixing its nodes in positions roughly corresponding to the players' formation on the pitch (see Figure \ref{fig:spvsned}). 

The passing network is by all means an oversimplification of a football game, as players do not remain in static positions during games. However, the network, with its arrows represented in various thickness and hue, does provides an immediate insight into a teams' tactics. It can be used, for example, to determine areas of the pitch that are favored or neglected, whether the team tends to use or abuse short distance or long distance passes, and whether a player is not intervening enough in a game. The network can also be used by a team to detect under-performing players, fix weak spots, detect potential problems between teammates who are not passing the ball as often as their position dictates, as well as to detect weaknesses in rivals.

This basic visual analysis can be made more quantitative by computing \textit{global} network invariants, which characterize a team as a whole, or \textit{local} invariants, which provide insight about individual players. The computation of some of these invariants is described in the next section. They all rely on the \emph{(weighted) adjacency matrix} $A$, having as entry $A_{ij}$ the number of passes from player $i$ to player $j$.

The weight (number of passes) will be used as a measure of the strength of an arrow in the network and also to define a notion of \textit{distance} between players. This distance $d_{ij}$ is defined precisely as the \emph{geodesic distance} given by the length of the shortest path connecting the nodes $i$ and $j$, where the length of a path is obtained by adding the lengths $l_{ij}$ of the arrows according to
\[
	l_{ij} = 
        \begin{cases}
		0 & \text{if $i=j$} \\
		\frac{1}{A_{ij}} & \text{if $i \neq j$}.
	\end{cases}
\]
The length of an arrow between two players is considered infinite if they do not pass the ball to each other. It is worth noting that our definition of distance does not need to be symmetric (i.e., one can have $d_{ij} \neq d_{ji}$), and is not necessarily correlated with the physical distance between the players in the field. For some computations we will also use the non-weighted adjacency matrix $\mathcal{E} = (\varepsilon_{ij})$, where
\[
\varepsilon_{ij} = \begin{cases}
  1 & \text{if $A_{ij}\neq 0$} \\
  0 & \text{if $A_{ij}= 0$}.
  \end{cases}
\]

An example of global invariant is the \emph{node connectivity}, defined as the minimum number of nodes one needs to remove in a network to make it disconnected. By the very nature of football games, this is not a good invariant to consider because the passing network of a team is usually very close to be complete, and thus has a very high degree of node connectivity. More useful is the \emph{edge connectivity}, defined as the minimum number of edges one needs to remove to make the network disconnected. This gives us a good measure of the game-play robustness, as it represents the smallest number of passes that need to be intercepted to interrupt a team's `natural flow' and to isolate a subset of its players, either by preventing the ball from reaching them or, if the edge connectivity is computed without the passes' directions, by completely isolating them from the rest of their teammates.

\section{Player performance}

The individual contribution of a player in a team can be inferred from local network invariants of the passing network and, in particular, from centrality measures, which define the \emph{relevance} or \emph{popularity} of a player according to different parameters. We define in this section three of these measures and discuss their meaning in the context of football.

\subsection{Closeness}
\label{sec:closeness-centrality}

The \textit{closeness centrality} or \textit{closeness score} of a player $i$ is one of the simplest notion of node centrality, defined as the inverse of the average geodesic distance of that node in the network \cite{Newman2010}: 
\[
C_i = \frac{20}{\sum_{j\neq i} d_{ij} + \sum_{j\neq i} d_{ji}}.
\]
For simplicity, we are giving equal weight to outgoing and incoming passes in this measure, but this can be adjusted for by throwing in arbitrary weights into the equation:
\[
C'_i = \frac{10}{w\sum_{j\neq i} d_{ij} + (1-w)\sum_{j\neq i} d_{ji}}.
\]
The closeness score provides a direct measurement on how easy it is to reach a particular player within a team. A high closeness score corresponds to a small average distance, indicating a well-connected player within the team.

 \subsection{Betweenness}
 \label{sec:betweenness}

A very different notion of centrality is \emph{betweenness centrality}, which measures the extent to which a node lies on paths between other nodes \cite{Newman2010}. This quantity is defined as the percentage of shortest paths that go through player $i$:
\[
C_B(i) = \frac{1}{90} \sum_{j\neq k\neq i} \frac{n^i_{jk}}{g_{jk}}
\]
where $n_{jk}^i$ is the number of geodesic paths from $j$ to $k$ going through $i$ and $g_{jk}$ is the total number of geodesic paths. The normalization factor $1/90$ ensures $0 \leq C_B(i) \leq 1$.

Betweenness does not measure how well-connected a player is, but rather how the ball-flow between \emph{other} players depends on that particular player $i$. It thus provides a measure of the impact of removing that player from the game, either by getting a red card or by being isolated by the rival's defense. A betweenness score of 0 means, in particular, that a player is not getting involved in the game, and so can be removed without much effect. 

From a tactical point of view, a team should seek betweenness scores that are evenly distributed among players: concentrated betweenness scores that are on the high side indicate a high dependence on few, too important players, whereas well distributed, low betweenness scores are an indication of a well-balanced passing strategy. 

Performance indicators based on betweenness centrality have been previously employed in the context of a team's activity \cite{Duch2010}. Further details on the computation of betweenness in directed networks with weights can be found in \cite{Opsahl2010}.

\subsection{Pagerank}
\label{sec:pagerank}
\textit{Pagerank centrality}, introduced in \cite{Brin1998}, is a recursive notion of `popularity' or importance which follows the principle that `a player is popular if he gets passes from other popular players'. Mathematically, pagerank centrality is defined by
\[
x_i = p \sum_{j\neq i} \frac{ A_{ji}}{L^{\textrm{out}}_{j}}x_j + q,
\]
where $L^{\textrm{out}}_j = \sum_{k} A_{jk}$ is the total number of passes made by player $j$, $p$ is a heuristic parameter representing the probability that a player will decide to give the ball away rather than keep it and go for a shot himself, and $q$ is a parameter awarding a `free' popularity to each player. Note that the pagerank score of a player depends on the scores of all his teammates. As a result, all pagerank scores in a team must be computed at the same time.

Pagerank centrality roughly assigns to each player the probability that he will have the ball after a reasonable number of passes has been made. If additional precision is required for this measurement, the probability $p$ can be replaced by player-dependent probabilities $p_i$, which would make more sense if certain players are more prone to keep the ball than others. In either case, the value of $p$ (or the $p_i$'s) does not come from the network alone, as it might in general be very different from one team to another, and should be determined by heuristics. As a proof of concept, in our analysis we will use a uniform value of $p = 0.85$ and $q = 1$ for all the teams studied.

\section{Clustering and communities}

An interesting aspect of football is how tightly players interact in a team. The notion of \emph{clustering} tells us precisely that: it is a measure of the degree to which nodes in a network tend to cluster together. 

The clustering coefficient of a node in a weighted network was originally defined in \cite{Barrat2004}; in our analysis we use a slight modification of that notion (see \cite{Saramaki2007} for a comparison of the different definitions) given by
\[
c_i^w = \frac{1}{u_i(u_i -1)} \sum_{j,k} \frac{\sqrt[3]{A_{ij}A_{kj}A_{ki}}}{\max(A)},
\]
where $u_i = \sum_j \varepsilon_{ij}$ is the number of passes made by player $i$, also known as the \emph{vertex out-degree}. 

The clustering coefficient accounts, technically speaking, for the \emph{transitivity} of the network by counting the percentage of all possible triangles containing the node $i$. To make this more precise, imagine that player $j$ wants to pass the ball to player $k$, but since this passing line is well defended, he has to go through player $i$ to reach player $k$, thus making the path $j\to i \to k$. If this is easy for them to do, i.e., if there is a high number of passes following the path $j\to i \to k$, then this translates in a high clustering score for player $i$, the one acting as middle-man. Conversely, if there is a large unbalance between the amount of passes involved in team, then the clustering coefficient will be lower. The average of all these coefficients gives us the \emph{global average clustering coefficient} for the team.

In addition to studying how clustered or fragmented a team is, we can compute the size of its \emph{maximal clique}, where a clique is a sub-network in which all the nodes are linked by an arrow. A clique in a team represents a subset of players that are all pairwise-connected by direct passes. A well connected team will present a very large maximal clique, meaning that almost everybody gets to pass the ball to everybody else, whereas the size will be smaller for more fragmented teams. The analysis of cliques is the basis for finding \emph{communities} within networks. 

Our initial attempt at studying communities within football teams has not provided any useful information, as the high degree of connectivity paired with the small number of nodes usually results in the existence of a single community that includes every player. This said, one might be able to use some suitable variation of the notion of community to overcome this inconvenience. We will postpone the study of this problem for future work.

\section{Results and analysis}

We present in this section the results of the computation of the different measures presented in the previous sections for the teams that have participated in the 2010 FIFA World Cup. For reasons of space, we only analyze the teams that have made it to the knock-out phase, and focus especially on the final qualifiers, Spain vs the Netherlands, and the third place qualifiers, Germany vs Uruguay. The passing networks and analysis of other teams can be found on our website\footnote{http://www.maths.qmul.ac.uk/\textasciitilde ht/footballgraphs/}.

\subsection{Data}

The passing data for the 2010 FIFA World Cup games was downloaded from the official FIFA website using a custom Python script. The passing networks were then constructed and analyzed using Sage \cite{Sage} and NetworkX \cite{Hagberg2008}. Graphics were created with Wolfram Research's Mathematica.

As FIFA only provides the aggregate data over all the games, the passing networks were computed by dividing the number of passes by the total number of games played by each team. This introduces artifacts in some cases. This issue can be taken care of by conducting a per-game analysis, which was unfortunately not possible in our case. 

\subsection{Teams in the last 16}

The centrality and clustering scores of the teams that made it to the last 16 stage are shown in Table \ref{table:top16}. Note that the betweenness and clustering scores are expressed as percentage of the theoretical maximum. 

The main point to note about these results is that Spain, the tournament winner and the team that arguably played the best football, has the highest number of passes, clustering and size of clique. It also has a high-end edge connectivity, while keeping a low betweenness score. All of this is a reflection of the `total football' or `\emph{tiki-taka}' style of Spain, in which well-connected players constantly pass the ball around. This is also confirmed by the passing network (see Figure \ref{fig:spvsned}) and the individual players' scores, discussed in the next subsection.

Other teams obtaining scores similar to Spain include the Netherlands (qualifying second in the tournament) and Brazil, followed by Argentina. At the lower end, Paraguay, with its low degree connectivity and and high betweenness, appears as a disconnected  team relying too much on a few players.

\begin{table}[t]
\centering
  \begin{tabular}{lcccccc}
    \textbf{Team} & $P$ & $k$ & $k_u$ & $\overline{c^w}$ & $\overline{C_B}$ & $Cq$ \\  \toprule
    Argentina & 227 & 4 & 5 & 27.9 & 2.7 & 8 \\  
    Brazil & \textbf{321} & \textbf{5} & \textbf{7} & 26.2 & 2.0 & 8 \\  
    Chile & 120 & 0 & 1 & 18.9 & \textbf{5.1} & 6 \\  
    England & 239 & 2 & 3 & 28.0 & 3.6 & 7 \\  
    Germany & 220 & 2 & 2 & 24.7 & 4.6 & 6 \\  
    Ghana & 184 & 3 & 4 & 15.5 & 3.5 & 8 \\  
    Japan & 180 & 1 & 5 & 28.9 & 3.3 & 8 \\  
    Korea Rep. & 227 & 3 & 5 & 24.4 & 2.6 & 8 \\  
    Mexico & 225 & 0 & 0 & 27.2 & 1.9 & 7 \\  
    Netherlands & 266 & \textbf{5} & \textbf{7} & \textbf{29.7} & 1.9 & 8 \\  
    Paraguay & 103 & 0 & 2 & 20.4 & \textbf{7.5} & 5 \\  
    Portugal & 175 & 3 & 4 & 14.6 & 4.1 & 7 \\  
    Slovakia & 166 & 3 & 6 & 18.5 & 3.0 & 7 \\  
    Spain & \textbf{417} & 3 & 5 & \textbf{30.0} & 1.9 & \textbf{9} \\  
    USA & 160 & 1 & 4 & 16.0 & 4.6 & 7 \\  
    Uruguay & 117 & 2 & 3 & 14.3 & 4.8 & 6 \\ \bottomrule
  \end{tabular}
  \caption{Data for the teams in the round of 16. $P$: average number of passes; $k$: edge connectivity; $k_u$: undirected connectivity; $\overline{c^w}$: average clustering; $\overline{C_B}$: average betweenness; $Cq$: largest clique. The highest two values (except for clique) are highlighted. }
    \label{table:top16}
\end{table}

\subsection{Spain vs the Netherlands}

Tables \ref{table:spain} to \ref{table:netherlands} show the closeness, betweenness, pagerank and clustering scores of the players of Spain and the Netherlands, respectively, in their formations used in the final.

Although there are some data artifacts due to the averaging of the data over several games (which again, would be sorted out by performing a per-game analysis), the overall conclusion that we reach from these results is that there is a high correlation between high scores in closeness, pagerank and clustering, which tend to confirm the general perception of the players' performance reported in the media at the time of the tournament. A remarkable example of this correlation is the high scores displayed by Xavi, arguably the leading player of the Spanish team.

On the Spanish side, one should also note that the betweenness scores are low and uniformly distributed -- a sign of a well-balanced passing strategy -- and consistently high clustering scores, showing that Spain is an extremely well-connected team, in which almost all players help each other by offering themselves as passing options. An exception is Pedro, whose low scores are  explained by the fact that he is a forward and was normally not playing games for their entire duration. Incidentally, note that forwards can almost always be identified as those players having the lowest closeness, betweenness and pagerank, as they are isolated players waiting to receive passes, as well as players who get replaced more often.

The scores of the Dutch team are, overall, close to those of the Spanish team, particularly the clustering, but there are some notable differences. First, there is a clear difference in the density of passes, seen in Figure \ref{fig:spvsned}. Second, the Dutch players are not as close to each other (as measured by $C_i$) and have pageranks that are more evenly distributed, thus showing that none has a predominant role in the passing scheme. Finally, Figure \ref{fig:spvsned} shows an unbalanced use of the pitch, giving a clear preference to the left side.

\begin{table}[t]
\centering
  \begin{tabular}{lcccc}
    \textbf{Player} & $C_i$ & $C_B(i)$ & $x_i$ & $c_i^w$ \\ \toprule 
    Casillas & 16.52 & 0.00 & 3.29 & 20.46 \\
    Pique & 17.32 & \textbf{3.92} & 11.46 & 30.70 \\
    Puyol & 16.32 & 2.86 & 7.92 & 27.12 \\
    Iniesta & 14.60 & 0.50 & 8.54 & 31.03 \\
    Villa & 8.68 & 0.50 & 5.89 & 23.96 \\
    Xavi & \textbf{18.28} & 1.19 & \textbf{14.66} & \textbf{46.47} \\
    Capdevila & 16.54 & \textbf{6.12} & 10.56 & 29.91 \\
    Alonso & 17.11 & 1.19 & 12.31 & \textbf{41.69} \\
    Ramos & 16.45 & 2.41 & 9.02 & 27.05 \\
    Busquets & \textbf{18.55} & 2.41 & \textbf{12.99} & 35.32 \\
    Pedro & 3.42 & 0.00 & 3.35 & 16.75 \\
    \bottomrule
  \end{tabular}
\caption{Player scores for Spain. The two highest scores are highlighted.}
\label{table:spain}
\end{table}

\begin{table}[t]
\centering
  \begin{tabular}{lcccc}
    \textbf{Player} & $C_i$ & $C_B(i)$ & $x_i$ & $c_i^w$ \\ \toprule 
    Stekelenburg & \textbf{16.34} & 0.32 & 7.63 & 28.35 \\
    Van Der Wiel & 14.43 & \textbf{2.97} & 9.79 & 31.39 \\
    Heitinga & 16.23 & 2.67 & \textbf{11.06} & 31.34 \\
    Mathijsen & \textbf{17.30} & 1.30 & 10.84 & 33.22 \\
    V. Bronckhorst & 15.74 & 1.12 & 10.07 & \textbf{37.00} \\
    Van Bommel & 12.46 & \textbf{3.08} & \textbf{11.19} & 32.36 \\
    Kuyt & 7.97 & 1.67 & 9.02 & 27.06 \\
    De Jong & 10.95 & 2.73 & 9.28 & 28.36 \\
    Van Persie & 6.89 & 2.92 & 5.88 & 20.13 \\
    Sneijder & 10.91 & 2.17 & 10.32 & \textbf{33.77} \\ 
    Robben & 5.91 & 0.16 & 4.91 & 23.91 \\
    \bottomrule
  \end{tabular}
\caption{Player scores for the Netherlands.}
\label{table:netherlands}
\end{table}

\subsection{Germany vs Uruguay}

Tables \ref{table:germany} to \ref{table:uruguay} shows the same data as in the previous subsection but now for Germany and Uruguay. The formations used in this case are those of the semifinals. 

The results for these two teams point to a major difference in their connectedness: Germany, with its closeness scores and high clustering, is overall more connected than Uruguay. However, as its pagerank scores are not as spread out as those of Uruguay, it seems to be depending more on the efforts of a few players to pass the ball around. Lahm and Schweinsteiger, in particular, are playing a central role in their team, not dissimilar to Xavi's.

\begin{table}[t]
\centering
  \begin{tabular}{lcccc}
    \textbf{Player} & $C_i$ & $C_B(i)$ & $x_i$ & $c_i^w$ \\ \toprule 
    Neuer & 7.58 & 0.37 & 4.74 & 21.54 \\
    Friedrich & 9.29 & 3.55 & 10.08 & 24.99 \\
    Khedira & 8.70 & 10.58 & 11.38 & 26.31 \\
    Schweinsteiger & 10.28 & \textbf{13.17} & \textbf{17.32} & 27.35 \\
    \"Ozil & 7.54 & 4.34 & 10.05 & 22.62 \\
    Podolski & 4.91 & 0.22 & 6.66 & \textbf{30.21} \\
    Klose & 0.92 & 0.00 & 2.48 & 14.34 \\
    Trochowski & 3.00 & 0.00 & 2.85 & \textbf{33.02} \\
    Lahm & \textbf{10.60} & \textbf{11.83} & \textbf{14.65} & 24.56 \\
    Mertesacker & \textbf{10.81} & 3.42 & 13.27 & 26.71 \\
    Boateng & 6.85 & 3.63 & 6.52 & 19.85 \\
    \bottomrule
  \end{tabular}
\caption{Player scores for Germany.}
\label{table:germany}
\end{table}

\begin{table}[t]
\centering
  \begin{tabular}{lcccc}
    \textbf{Player} & $C_i$ & $C_B(i)$ & $x_i$ & $c_i^w$ \\ \toprule 
    Muslera & 0.88 & 1.98 & 4.62 & 9.96 \\
    Godin & 1.80 & 4.20 & 8.37 & 11.99 \\
    Gargano & 0.76 & 0.37 & 2.98 & 7.33 \\
    Victorino & 1.75 & 0.88 & 5.59 & 14.28 \\
    Cavani & 1.61 & 10.22 & 10.22 & 13.68 \\
    Forlan & 2.08 & \textbf{10.29} & \textbf{13.12} & 15.02 \\
    A. Pereira & 1.90 & 3.75 & 9.01 & 16.45 \\
    Perez & \textbf{2.36} & \textbf{10.63} & \textbf{15.25} & 19.12 \\
    M. Pereira & 2.28 & 1.51 & 12.21 & \textbf{20.07} \\
    Arevalo & \textbf{2.45} & 5.85 & \textbf{13.12} & \textbf{19.83} \\
    Caceres & 1.34 & 3.65 & 5.52 & 9.51 \\
    \bottomrule
  \end{tabular}
\caption{Player scores for Uruguay.}
\label{table:uruguay}
\end{table}

\section{Further work}

The passing networks that we have presented provide an attractive visual summary or `snapshot' of a football team's style. The obvious limitation of these networks is of course that they are static. But, as we have seen, they can be complemented with the computation of centrality measures that provide useful information about the importance and connectedness of individual players, which might benefit coaches, sports journalists and their readers. 

There are many additional features that could be added to the networks to obtain a more detailed analysis. An immediate one would be to add an extra node representing the opponent's goal and consider shots instead of passes for arrows directed at the goal. This concept, with one node for shots in target and one for wide shots, has been previously used in \cite{Duch2010}. 

Another interesting aspect to consider would be to study the accuracy of passes by adding to each player a weight taking into account the probability for a pass coming from that player to be successful. There are different levels of complexity that one might want to get into here, as not all passes are equally likely to succeed or fail. But, as a first approximation, one might just want to use the percentage of completed passes as a measurement of accuracy.

Finally, let us mention that the defensive strength of a team could also be incorporated in the model by tracking passing interceptions and recovered balls.

\section*{Acknowledgements}
We thank J.J. Merelo G\"uerv\'os for introducing us to the possibility of using network theory in the context of football, C. Clanet for inviting us to the Euromech \emph{Physics of Sports} Conference and encouraging us to write down the results of our analysis as a paper, and L. Mahadevan for suggesting to study communities within football teams. We also thank FIFA for making the passing data publicly available during the duration of the World Cup and afterwards. The work of J. L\'opez Pe\~na was supported by an EU Marie-Curie fellowship (PIEF-GA-2008-221519) and a Spanish MCIM grant (MTM2010-20940-C02-01).


\nocite{Brillinger2007} \nocite{Duch2010} \nocite{Hughes2005}
\nocite{Joseph2006} \nocite{Lewis2009} \nocite{Luzuriaga2010}
\nocite{Newman2010} \nocite{Opsahl2010} \nocite{Yamamoto2011}
\nocite{Sage} \nocite{Hagberg2008} \nocite{Brin1998}

\end{document}